\documentclass[12pt,a4paper]{article}
\usepackage{latexsym}
\usepackage{amsfonts}
\usepackage[all]{xy}
\usepackage{amssymb, mathrsfs, amsfonts, amsmath}
\usepackage{amsbsy}
\newtheorem{thm}{Theorem}[section]
 \newtheorem{cor}[thm]{Corollary}
 \newtheorem{lem}[thm]{Lemma}

 \newtheorem{defn}[thm]{Definition}
 \newtheorem{rem}[thm]{Remark}
 
 \numberwithin{equation}{section}
 \newcommand{\noin}{\noindent}
 \newcommand{\hess}{\textrm{Hess}}

 \newcommand{\vol}{{\rm vol}}
\newcommand{\inj}{{\rm inj}}
\newcommand{\diver}{{{\rm{div}}}}

\newcommand{\grad}{{{\rm{grad}\,}}}

\begin{document}

\title{ Mean time exit and isoperimetric inequalities \\for minimal submanifolds of  $N\times \mathbb{R}$ }
\author{G. Pacelli Bessa\thanks{This work was completed when the authors were visiting the  Abdus Salan
International Center for Theoretical Physics, ICTP }\and J.
F\'{a}bio Montenegro }
\date{\today}

\maketitle
\begin{abstract}
\noin  Based on Markvorsen and Palmer's work on mean time exit and
isoperimetric inequalities we  establish
 slightly better isoperimetric inequalities and mean time exit estimates for minimal submanifolds of $N\times\mathbb{R}$.
 We also prove isoperimetric inequalities  for submanifolds  of Hadamard spaces with tamed second fundamental form.
\vspace{.2cm}


\end{abstract}

\section{Introduction}\hspace{.5cm}  The study of minimal surfaces in  product spaces $N\times \mathbb{R}$,
 where $N$ is a complete surface,  started with  H. Rosenberg in
\cite{rosenberg} and it has shown to be  a rich and interesting
theory, yielding a wealth of examples
 and  results, \cite{elbert-rosenberg}, \cite{hauswirth},
\cite{hauswirth-rosenberg}, \cite{meeks-rosenberg1}, \cite{nelli-rosenberg1}, \cite{rosenberg-book}.
 It also lead to the  study of constant mean curvature  surfaces in product spaces,  \cite{abresch-Rosenberg},
\cite{bessa-costa}, \cite{bessa-montenegro2}, \cite{earp},
 \cite{hoffman-lira-Rosenberg},
\cite{nelli-rosenberg2}, \cite{nelli-rosenberg3}.
 The classical theory of minimal surfaces in $\mathbb{R}^{3}$ guides the search  in this new theory,
 that depending on the geometry of  $N$,  yields very different results from their counterparts in the classical theory.
 In this spirit, based on the ideas of Markvorsen-Palmer,  we  study isoperimetric inequalities for minimal submanifolds
 of $ N\times \mathbb{R}$,  where $N$ is a complete Riemannian $n$-manifold  with sectional
 curvature $K_{N}\leq b$.
Markvorsen and  Palmer in  \cite{markvorsen-palmer} and
 \cite{palmer-1} proved  isoperimetric inequalities for extrinsic geodesic
 balls of proper minimal submanifolds  of Riemannian manifolds  with sectional curvature bounded above. To be precise,  let
        $\varphi: M\hookrightarrow W$  be a proper minimal immersion of an   $m$-dimensional
     manifold $M$ into
      a Riemannian
     $n$-manifold $W$ with sectional curvature $K_{W}\leq b$ and let
       $B_{W}(R)$ be a
       geodesic ball of $W$ centered at a point $p=\varphi (q)$ with radius $R\leq \min \{\inj_{W}(p), \pi/2\sqrt{b}\}$, where
      $\pi/2\sqrt{b}=\infty$ if $b\leq 0$ and $\inj_{W}(p)$ is the injectivity radius at
      $p$. The
        extrinsic geodesic ball of radius $R$ centered at $p$, denoted by  $D(R)$,  is defined to be the connected
       component of $\varphi (M)\cap B_{W}(R)$ containing  $p$. The isoperimetric inequalities proved by Markvorsen and Palmer are the following.
  \begin{thm}[Markvorsen-Palmer, \cite{markvorsen-palmer}, \cite{palmer-1}]\label{palmer}
\begin{itemize}\item[]\item[i.] If $b\leq 0$ then
\begin{equation}\frac{\vol_{m-1}(\partial D(R))}{\vol_{m}(D(R))}\geq
\frac{\vol_{m-1}(\partial B_{\mathbb{N}^{m}(b)}(R))}{\vol_{m}( B_{\mathbb{N}^{m}(b)}(R))}\label{eqPalmer-1}
\end{equation}
\item[ii.] If $b>0$ then
\begin{equation}\label{eqPalmer-2}\frac{\vol_{m-1}(\partial D(R))}{\vol_{m}(D(R))}
\geq m\displaystyle\frac{C_{b}}{S_{b}}(R).
\end{equation}
\end{itemize}\noindent
Where $(C_{b}/S_{b})(R)$ is the (constant) mean curvature of the geodesic
sphere $\partial B_{\mathbb{N}^{m}(b)}(R)$ of radius $R$ in $\mathbb{N}^{m}(b)$
 and the functions $S_{b}$ and $C_{b}$ are defined in (\ref{eqSk}). Moreover,
 equality in item i. implies that $D(R)$  is a minimal cone in $W$.
\end{thm}We consider  a minimal immersion
  $\varphi :M \hookrightarrow N\times \mathbb{R}$  of $m$-dimensional
     manifold $M$ into the product space $N\times \mathbb{R}$, where $N$ is a complete Riemannian manifold with sectional curvature $K_{N}\leq b$.
     Let
 $K\subset\varphi (M)$ be a connected compact set and let $r_{_{K}}={\rm rad}(\pi_{1}(K))$ be the radius of the set $\pi_{1}(K)$, where
   $\pi_{1}:N\times \mathbb{R}\to N$ is the projection on the first factor. Denote by
    $p_{_{K}}\in N$  the barycenter of $\pi_{1}(K)$ and suppose that $r_{_{K}}< \min \{ \inj_{N}(p_{_{K}}), \pi/2\sqrt{b}\}$. We prove a slightly better isoperimetric inequality when $b<0$.
\begin{thm}\label{thmA}
 If $b\leq 0$ then
\begin{equation}\frac{\vol_{m-1}(\partial K)}{\vol_{m}(K)}\geq
\frac{\vol_{m-2}(\partial
B_{\mathbb{N}^{m-1}(b)}(r_{_{K}}))}{\vol_{m-1}(
B_{\mathbb{N}^{m-1}(b)}(r_{_{K}}))}\cdot\label{eqThm1-i}
\end{equation}

\end{thm}
\begin{rem}This  isoperimetric inequality is sharp if we consider arbitrary compact sets $K$.
Consider the totally geodesic
 embedding $\varphi :\mathbb{H}^{m-1}\times \mathbb{R}\hookrightarrow \mathbb{H}^{n}(-1)\times \mathbb{R}$ given by
  $\varphi (x,t)=(x,t)$ and a family of compact sets $K_{i}=B_{\mathbb{H}^{m-1}}(R)\times [-i,i]$, $i=1,2,\ldots$. We then have that $$\frac{\vol_{m-1}(\partial K_{i})}{\vol_{m}(K_{i})}=\frac{\vol_{m-2}(\partial B_{\mathbb{H}^{m-1}}(R) )}{\vol_{m-1}(B_{\mathbb{H}^{m-1}}(R))}+\frac{1}{i} \rightarrow \frac{\vol_{m-2}(\partial B_{\mathbb{H}^{m-1}}(R) )}{\vol_{m-1}(B_{\mathbb{H}^{m-1}}(R))}\cdot$$
On the other hand, if we consider extrinsic geodesic balls $D(R)$ of $M\hookrightarrow  N\times \mathbb{R}$, $K_{N}\leq -1$, then Markvorsen-Palmer's estimate (\ref{eqPalmer-1}) is better if
  $m=2$. Since
 $$\displaystyle\frac{\vol_{1}(\partial D(R))}{\vol_{2}(D(R))}
 \geq
  \frac{\vol_{1}(\partial B_{\mathbb{R}^{2}}(R))}{\vol_{2}(B_{\mathbb{R}^{2}}(R))}=
 \frac{2}{R}>\frac{1}{R}=\frac{\vol_{0}(\partial
 B_{\mathbb{H}^{1}}(R))}{\vol_{1}(B_{\mathbb{H}^{1}}(R))}\cdot$$
For $m\geq 3$ there exists $R_{m}$ such that if  $R\geq R_{m}$ then we have   $$\displaystyle\frac{\vol_{m-1}(\partial
D(R))}{\vol_{m}(D(R))} \geq
   \frac{\vol_{m-2}\partial B_{\mathbb{H}^{m-1}(-1)}(R)}{\vol_{m-1} B_{\mathbb{H}^{m-1}(-1)}(R)}=
   \frac{\sinh(R)^{m-2}}{\int_{0}^{R}\sinh(s)^{m-2}ds}\geq \frac{m}{R}=  \frac{\vol_{m-1}(\partial
B_{\mathbb{R}^{m}}(R))}{\vol_{m}(B_{\mathbb{R}^{m}}(R))}\cdot$$ In fact, a rough estimate gives  $$ \frac{\sinh(R)^{m-2}}{\int_{0}^{R}\sinh(s)^{m-2}ds}\geq (m-2)\cdot    \frac{(e^{R}-1)^{m-2}}{e^{(m-2)R}-1}$$
Just let $R_{m}$ be such that $$ (m-2)\cdot  \frac{(e^{R_{m}}-1)^{m-2}}{e^{(m-2)R_{m}}-1}=\frac{m}{R_{m}}\cdot$$
\end{rem}

   Our next result gives upper bounds for the   isoperimetric quotients for extrinsic geodesic balls of submanifolds with tamed second fundamental form in Hadamard spaces with bounded sectional curvature.

   Let  $\varphi: M
\hookrightarrow N$ be  an isometric immersion of a complete Riemannian
$m$-manifold $M$ into a Hadamard $n$-manifold $N$  with
sectional curvature bounded above $K_{N}\leq b \leq 0$. Fix a
point $x_0 \in M$ and let $\rho_{M} (x) = {\rm dist}_{M}(x_0, x)$ be the
 distance function on $M$ to $x_0$.

 \vspace{2mm}
 \noindent Let
$\{C_{i}\}_{i=1}^{\infty}$ be an exhaustion sequence of $M$ by
compacts sets with $x_0 \in C_1$ and define a non-increasing
sequence $a_1 \geq a_2 \geq \cdots \geq 0$ by
\begin{eqnarray*}
\begin{array}{ccl}
a_i =  \sup \left \{ \displaystyle \frac{S_{b}(\rho_{M} (x))}{
C_{b}(\rho_{M} (x))}\cdot  \Vert \alpha (x)\Vert, \, x \in M
\backslash C_i \right \},
\end{array}
\end{eqnarray*}
where
\begin{equation}\label{eqSk}
 S_{b} (t)=\left \{
\begin{array}{ccl}
\displaystyle \frac {1}{\sqrt{-b}}\sinh(\sqrt{ -b}\,t),&if& b<0 \\
t, &if& b=0, \\
\displaystyle \frac {1}{\sqrt{b}}\sin(\sqrt{ b}\,t),&if& b>0 \, \, {\rm and }\,\, t<\pi/2\sqrt{b}
\end{array} \right.
\end{equation}
$ C_{b}(t)= S_{b}'(t)$ and $\alpha(x)$ is the  second
fundamental form of $\varphi(M)$ at $\varphi(x)$. It is clear that
the limit  $a(M)=\lim_{i\to \infty}a_{i}\in [0, \infty]$ does not
depend  on the exhaustion sequence $\{C_{i}\}_{i=1}^{\infty}$ nor on the base point $x_{0}$.

\vspace{2mm}
\begin{defn}\label{def1} An immersion $\varphi: M \hookrightarrow N$  of a complete Riemannian $m$-manifold
$M$ into a Hadamard $n$-manifold $N$ with sectional curvature
 $K_{N}\leq b\leq 0$ has tamed second fundamental form if $a(M)< 1$.\end{defn}
 \vspace{2mm}

 \noindent  Submanifolds of $\mathbb{R}^{n}$ with tamed second fundamental form were studied by the authors and L. Jorge  in \cite{bessa-jorge-montenegro} where we showed that complete submanifolds with
  tamed fundamental form of the $\mathbb{R}^{n}$ are proper and has finite topology.  Silvana Costa \cite{costa}
   extended this result  to submanifolds  of Hadamard manifolds with tamed second fundamental form. Here we give upper bounds for the isoperimetric quotients of extrinsic geodesic balls in Hadamard manifolds. We  prove the following theorem.

\begin{thm}\label{thm2} Let $\varphi: M \hookrightarrow N$  be a complete immersed $m$-submanifold $M$ with
 tamed second fundamental form of an $n$-dimensional Hadamard manifold $N$ with bounded sectional curvature
 $b_{1}\leq K_{N}\leq b_{2}\leq 0$. For a  given $c\in (a(M),1)$ there exists positive  constants $r_{0}=r_{0}(b_{2}, c)$, $ B=B(b_{2},c)<1$
 such that for extrinsic geodesic balls $D(R)$ with radius $R\geq r_{0}$ we have
 \begin{equation} \frac{{\rm vol}_{m-1}(\partial D(R))}{{\rm vol}_{m}(D(R))}\leq \frac{1+ \sqrt{-b_{1}}\cdot R\cdot \coth(\sqrt{-b_{1}}\cdot R)+ \Lambda }{R\cdot \sqrt{1-B^{2}}}
 \end{equation}
Where $\Lambda$ is a constant depending on  $c$, $r_{0}$, $R$ ,$b_{2}$ and $\sup_{B_{N}(r_{0})}\vert H\vert$.
\end{thm}

  \section{Proof of the results}
\subsection{Basic formulas}

Let $\varphi : M \hookrightarrow W$ be an isometric immersion
$M$ and $W$ are  Riemannian
 manifolds. Consider a smooth function $g:W \rightarrow \mathbb{R}$ and the composition $f=g\,\circ\,
 \varphi :M \rightarrow \mathbb{R}$.
Let $\nabla
$ and $\overline{\nabla}$ be the Riemannian connections on $M$ and
$W$ respectively,  $\alpha (q) (X,Y) $ and $\hess\,f(q)\,(X,Y)$  be
respectively the second fundamental form of the immersion $\varphi $
and
 the Hessian of $f$ at $q\in M$,  $X,Y \in T_{p}M$.  Identifying $X$ with $d\varphi (X)$ we have at $q\in M$ and  for every
 $X\in T_{q}M$ that
 \begin{equation}\hess\,f (q)  \,(X,Y)= \hess\,g (\varphi (q))\,(X,Y) +
 \langle \grad\,g\,,\,\alpha (X,Y)\rangle_{\varphi (q)}.
\label{eqBF2}
\end{equation}
 Taking the trace in (\ref{eqBF2}), with respect to an orthonormal basis $\{ e_{1},\ldots e_{m}\}$
 for $T_{q}M$, we have  that
\begin{eqnarray}
\Delta \,f (q)                & = & \sum_{i=1}^{m}\hess\,f (q)  \,(e_{i},e_{i})\nonumber \\
                             & = & \sum_{i=1}^{m}\hess\,g (\varphi (q))\,(e_{i},e_{i}) + \langle \grad\,g\,,\,
                             \sum_{i=1}^{m}\alpha (e_{i},e_{i})\rangle.\label{eqBF3}
\end{eqnarray}
 The formulas (\ref{eqBF2}) and (\ref{eqBF3}) are  well known in the literature,
 see  \cite{jorge}.
Another important tool   is the Hessian Comparison Theorem, see \cite{schoen}.

\begin{thm}[Hessian Comparison Theorem] Let $W$ be a  complete Riemannian $n$-manifold  and
$y_{0},y_{1} \in W $.   Let
 $\gamma:[0,\,\rho_{W} (y_{1}) ]\rightarrow M$ be a minimizing geodesic joining $y_{0}$ and $y_{1}$ where $\rho_{W} $
 is the  distance function  to $y_{0}$ on $W$. Let $K_{\gamma}$ be the sectional curvatures of $W$ along $\gamma$ and let $c=\inf K_{\gamma}$ and $b=\sup K_{\gamma}$.
 The Hessian of  $\rho_{W}$ at $y=\gamma( \rho_{W}(y))$  for any $X\in T_{y}W$, $X\perp \gamma'(\rho_{W}
(y))$,  satisfies

\begin{equation}
\frac{C_{c}}{S_{c}}(\rho_{W}(y))\cdot\Vert X\Vert^{2}\geq  Hess\,\rho_{W}(y)(X,X)\geq  \frac{C_{b}}{S_{b}}(\rho_{N}(y))\cdot\Vert X\Vert^{2},\label{eqBF6}
\end{equation}whereas $Hess\,\rho_{W}(y) (\gamma ',\gamma ') =0$.
\label{thmHess}
\end{thm}
\subsection{Mean time exit from minimal submanifolds of $N\times \mathbb{R}$}
\noindent  Let $\varphi :M\hookrightarrow W$ be a complete, minimal,  properly  immersed $m$-submanifold of a
 complete Riemannian manifold $W$  with sectional curvature  $K_{W}\leq b$. Let $D(R)$ be an extrinsic geodesic
 ball centered at $p=\varphi (q)$ with radius $R$ and $\rho_{_{W}}(x)={\rm dist}_{W}(p,x)$. Let $E(x)$ be the mean time of the first exit from $D(R)$  of a
  particle in Brownian motion starting at $x\in D(R)$ and denote by
   $E^{m}_{b}(\tilde{x})=E^{m}_{b}(\vert \tilde{x}\vert)$  the mean time of the first exit from
   $B_{\mathbb{N}^{m}(b)}(R)$ of a particle in Brownian motion starting at $\tilde{x}\in B_{\mathbb{N}^{m}(b)}(R)$,
   $\vert \tilde{x}\vert ={\rm dist}_{\mathbb{N}^{m}(b)}(0,\tilde{x})$.  Markvorsen,  \cite{markvorsen2} proved the following theorem.
  \begin{thm}[Markvorsen's mean time exit comparison theorem]\label{mark}
   \begin{itemize}\item[]
   \item[i.] If the sectional curvature $b\geq K_{W}\geq \kappa\geq 0$
    then $E(x)\geq E^{m}_{\kappa}(\rho_{_{W}}(x))$.
 \item[ii.] If If the sectional curvature $0\geq b\geq K_{W}$ then $E(x)\leq E^{m}_{b}(\rho_{_{W}}(x))$. \end{itemize}
  \end{thm}We have a version of Markvorsen's mean time exit comparison theorem for compact sets of minimal submanifolds of $N\times \mathbb{R}$.
Let $K\subset\varphi (M)$ be compact set in a minimal $m$-submanifold of  $ N\times\mathbb{R}$,  where $N$
is a Riemannian $n$-manifold  with sectional curvature $K_{N}\leq b$. Let $r_{_{K}}$ and $p_{_{K}}$ be respectively
 the radius  and  barycenter of $\pi_{1}(K)$. Suppose that $r_{_{K}}< \min\{ \inj_{N}(p_{_{K}}), \pi/2\sqrt{b}\}$.
  Denote by $E(x)$ the mean time of the first exit   from $K$ of a particle in Brownian motion starting at
   $x\in K$ and  by $E_{b}(\tilde{x})=E^{m-1}_{b}(\vert \tilde{x}\vert)$  the mean time of the first
    exit from $B_{\mathbb{N}^{m-1}(b)}(r_{_{K}})$ of a particle in Brownian motion starting at
     $\tilde{x}\in B_{\mathbb{N}^{m-1}(b)}(r_{_{K}})$, $\vert \tilde{x}\vert ={\rm dist}_{\mathbb{N}^{m-1}(b)}(0,\tilde{x})$.

  \noindent
  We prove the following comparison theorem.

   \begin{thm}\label{thmMTE}\begin{itemize}\item[]
    \item[i.]If $K_{N}\leq b\leq 0$ then $E(x)\leq E_{b}(\rho_{N}(\pi_{1}(x))).$\item[ii] If $K_{N}\geq \kappa \geq 0$, suppose that the immersion $\varphi $ is proper and $K=\varphi(M)\cap (B_{N}(R)\times \mathbb{R})$  is a compact set, then $ E(x)\geq E_{\kappa}(\rho_{N}(\pi_{1}(x)))$. Where $\rho_{N}(\pi_{1}(x))={\rm dist}_{N}(p_{_{K}},\pi_{1}(x))$.
\end{itemize}
  \end{thm}
\begin{rem}The statement of this theorem is somewhat surprising. For instance,  consider the totally geodesic
 embedding $\varphi :\mathbb{H}^{m-1}\times \mathbb{R}\hookrightarrow \mathbb{H}^{n}(-1)\times \mathbb{R}$ given by
  $\varphi (x,t)=(x,t)$ and $K=B_{\mathbb{H}^{m-1}}(R)\times [-L,L]$. It does not matter how large is $L$, the mean time exit of $K$ can not exceed $E_{-1}(R)$. The particle in Brownian motion can not move upward for too long. It is drifted  horizontally to the boundary.
\end{rem}

\subsection{Proof of Theorem \ref{thmMTE}}

 We denoted by $E(x)$ the mean time of the first
exit
 from $K$ of a particle in Brownian motion starting at $x$ and by $E_{b}(\tilde{x})$  the mean time from the
 first exit of the geodesic ball $B_{\mathbb{N}^{m-1}(b)}(r_{_{K}})$ of a particle in Brownian motion starting at
  $\tilde{x}$.
 A remark from Dynkin \cite{dynkin}{ vol 2, p.51}  states that the functions $E$ and $E_{b}$ satisfies the Dirichlet
  boundary problem.
  \begin{equation} \label{eqDynkin-A}\left\{\begin{array}{rcrll}\triangle_{K}E&=&-1 &{\rm  in} & K\\
                                             E & =&0 &{\rm on}& \partial K\end{array}\right. \,\& \,\,
                                              \left\{\begin{array}{rcrll} \triangle_{ \mathbb{N}^{m-1}(b)}E_{b}
                                              &=&-1 &{\rm  in} & B_{\mathbb{N}^{m-1}(b)}(r_{_{K}})\\
                                             E_{b} & =&0 &{\rm on}& \partial  B_{\mathbb{N}^{m-1}(b)}(r_{_{K}})
                                             \end{array}\right.\end{equation}It is known that $E_{b}$
                                              is a radial function $E_{b}(\tilde{x})=
                                              E_{b}(\vert \tilde{x}\vert)$,
                                              $\vert \tilde{x}\vert={\rm dist}_{\mathbb{N}^{m-1}(b)}(0, \tilde{x})$.
                                              Let $\bar{E}_{b}$  be the transplant of $E_{b}$
                                               to $B_{N}(r_{_{K}})\times \mathbb{R}$ defined by
                                                $\bar{E}_{b}(x)=E_{b}\circ \rho_{N}\circ \pi_{1}(x)$,
                                                where $\pi_{1}:N\times \mathbb{R}\to N$ is the projection
                                                on the first factor. We have that $\bar{E}_{b}\vert_{K}=
                                                \bar{E}_{b}\circ \varphi$.
 Following Markvorsen \cite{markvorsen2} we define $F_{b}:[0,\infty)\to [0, \infty) $ by
  \begin{equation} F_{b}(t)=\left\{\begin{array}{lll} \displaystyle \frac{1}{b}(1-\cos(\sqrt{b}\cdot t) &if & b>0\\
  && \\\displaystyle \frac{t^{2}}{2} & if & b=0\\&&\\
  \displaystyle \frac{1}{b}(1-\cosh(\sqrt{-b}\cdot t) &if & b<0 \end{array}\right.
  \end{equation}Observe that $F_{b}$ satisfies $ F_{b}''(t)-(C_{b}/S_{b})(t)F_{b}'(t)=0$ for all $t\geq 0$.

  \vspace{2mm}
  \noindent Let $s=F_{b}(\rho_{N}\circ \pi_{1})$ and define $\pmb{ \bar{E}}_{b}(s)$ by $\pmb{ \bar{E}}_{b}(s(x))=\bar{E}_{b}(x)$.   Computing $\triangle_{K}\bar{E}_{b}\circ \varphi$ at any point $x\in B_{N}(R)$   we obtain, (see \ref{eqBF3})

  \begin{eqnarray}\label{eq2.7}\triangle_{K}\bar{E}_{b}\circ \varphi(x)
  &=& \sum_{i=1}^{m}\hess_{(N\times \mathbb{R})}\bar{E}_{b}(y)(X_{i},X_{i})\nonumber \\
  &=&\sum_{i=1}^{m}\hess_{(N\times \mathbb{R})}\,E_{b}\circ\rho_{N}\circ \pi_{1}(y)(X_{i},X_{i})\nonumber \\
  &=& \sum_{i=1}^{m}\hess_{N}(E_{b}\circ\rho_{N}\circ \pi_{1})(y)(X_{i}, X_{i})\\
  &=& \sum_{i=1}^{m}\hess_{N}\pmb{\bar{E}}_{b}(s(y))(X_{i}, X_{i})\nonumber\\
  &=& \sum_{i=1}^{m}\left[ \pmb{\bar{E}}_{b}''(s(y))\langle \grad s, X_{i}\rangle^{2} + \pmb{\bar{E}}_{b}'(s(y))\hess_{N}\,s(y)(X_{i},X_{i})\right]\nonumber
  \end{eqnarray}Where $\{X_{i}\}$ is an orthonormal basis for $T_{y}\varphi (M)$, $y=\varphi (x)$. Let $\{\partial/\partial \rho_{N}, \partial/\partial \theta_{1}, \ldots,\partial/\partial \theta_{n}\}$ be an orthonormal basis for $T_{\pi_{1}(y)}N$ from polar coordinates and $\partial/\partial t$ is the a tangent to the $\mathbb{R}$ factor. We choose  $\{X_{i}\}$ in the following way.
  \begin{eqnarray}\label{basis} X_{i}&=&\alpha_{i}\cdot \partial/\partial \rho_{N}+\beta_{i}\cdot\partial/\partial t + \sum_{j=1}^{n-1}\gamma_{j}^{i} \cdot\partial/\partial \theta_{j}\end{eqnarray} \begin{eqnarray} \alpha_{i}^{2}+\beta_{i}^{2}+\sum_{j=1}^{n-1}(\gamma_{j}^{i})^{2}&=1&.\end{eqnarray}
We compute $ \sum_{i=1}^{m}\hess_{N} \,s (X_{i},X_{i})$ taking in account  the Hessian Comparison Theorem and the fact $F_{b}''(t)-F_{b}'(t)(C_{b}/S_{b})(t)=0$.
\begin{eqnarray}\label{hessiano-de-s-A}\sum_{i=1}^{m}\hess_{N} \,s (X_{i},X_{i})&=&\hspace{-1mm} F_{b}''(\rho_{N})\sum_{i=1}^{m}\langle \grad \rho_{N}, X_{i}\rangle^{2} +  F_{b}'(\rho_{N})\sum_{i=1}^{m}\hess \rho_{N} (X_{i},X_{i})\nonumber  \\
&=&F_{b}''(\rho_{N})\sum_{i=1}^{m}\alpha_{i}^{2} + F_{b}'(\rho_{N})\sum_{i=1}^{m}\sum_{j=1}^{n-1}(\gamma_{j}^{i})^{2}\hess_{N}\rho_{N}(\partial/\partial \theta_{j},\partial/\partial \theta_{j}) \\
 &\geq & F_{b}''(\rho_{N})\sum_{i=1}^{m}\alpha_{i}^{2} +  F_{b}'(\rho_{N})\frac{C_{b}}{S_{b}}(\rho_{N})\sum_{i=1}^{m}(1-\alpha_{i}^{2}-\beta_{i}^{2})\nonumber \\
&=& \left(F_{b}''(\rho_{N})-F_{b}'(\rho_{N})\frac{C_{b}}{S_{b}}(\rho_{N})\right)\sum_{i=1}^{m}\alpha_{i}^{2} +F_{b}'(\rho_{N})\frac{C_{b}}{S_{b}}(\rho_{N})(m-\sum_{i=1}^{m}\beta_{i}^{2})\nonumber \\
&\geq &(m-1)F_{b}'(\rho_{N})\frac{C_{b}}{S_{b}}(\rho_{N})\nonumber
\end{eqnarray} Thus \begin{eqnarray}\label{hessiano-de-s}
\sum_{i=1}^{m}\hess_{N} \,s(X_{i},X_{i})&\geq & (m-1)\cdot F_{b}'(\rho_{N})\frac{C_{b}}{S_{b}}(\rho_{N}).
\end{eqnarray}
Recall that the Laplacian of the canonical metric  $dt^{2}+ S_{b}^{2}(t)d\theta^{2}$ of the space form $\mathbb{N}^{m-1}(b)$ is given by $\triangle_{\mathbb{N}^{m-1}(b)}=\partial^{2}/\partial t^{2} + (m-2)(C_{b}/S_{b}) \partial/\partial t + (1/S_{b}^{2}(t))\triangle_{\mathbb{S}^{m-2}}$.

\noindent Therefore
$$\triangle_{\mathbb{N}^{m-1}(b)}s=\triangle_{\mathbb{N}^{m-1}(b)}F_{b}(\rho_{N})=
F_{b}''(\rho_{N})+(m-2)\frac{C_{b}}{S_{b}}F_{b}'(\rho_{N})=(m-1)\frac{C_{b}}{S_{b}}F_{b}'(\rho_{N}).$$
In \cite{markvorsen2}, Proposition 4, Markvorsen proved that $\pmb{\bar{E}}_{b}'(s)<0$ for all $b$ and $\pmb{\bar{E}}_{b}''(s)>0$ if $b<0$, $\pmb{\bar{E}}_{b}''(s)=0$ if $b=0$ and $\pmb{\bar{E}}_{b}''(s)<0$ if $b>0$.

 \noindent Therefore from (\ref{eq2.7}) we have
 \begin{eqnarray}\label{eq3.26}\triangle_{K}\bar{E}_{b}\circ \varphi(x)
  &=&  \sum_{i=1}^{m}\left[ \pmb{\bar{E}}_{b}''(s)\langle \grad s, X_{i}\rangle^{2} + \pmb{\bar{E}}_{b}'(s)\hess_{N}\,s(X_{i},X_{i})\right]\nonumber\\
  &\leq & \pmb{\bar{E}}_{b}''(s)\vert \grad_{N} s\vert^{2} +  \pmb{\bar{E}}_{b}'(s)\cdot (m-1)\cdot F_{b}'(\rho_{N})\displaystyle \frac{C_{b}}{S_{b}}(\rho_{N})\nonumber \\
 &&\\ &=&  \pmb{\bar{E}}_{b}''(s)\vert \grad_{\mathbb{N}^{m-1}(b)} s\vert^{2} +\pmb{\bar{E}}_{b}'(s)\triangle_{\mathbb{N}^{m-1}(b)}s\nonumber\\
 &&\nonumber\\
 &=& \triangle_{\mathbb{N}^{m-1}(b)}\pmb{\bar{E}}_{b}=-1=\triangle_{K}E\nonumber
 \end{eqnarray}
   Then $\triangle_{K}(\bar{E}_{b}-E)\leq 0$ with $(\bar{E}_{b}-E)\vert_{\partial K}=\bar{E}_{b}\vert_{\partial K}\geq 0$. Thus  $\bar{E}_{b}\geq E$ in $K$.
\vspace{2mm}

\noindent If $K_{N}\geq \kappa>0$ then \begin{eqnarray}\label{eq3.30}
\sum_{i=1}^{m}\hess_{N} \,s(X_{i},X_{i})&\leq & (m-1)\cdot F_{b}'(\rho_{N})\frac{C_{b}}{S_{b}}(\rho_{N}).
\end{eqnarray}and $\pmb{\bar{E}}_{\kappa}''(s)<0$. The same reasoning as before shows that $\triangle_{K}(\bar{E}_{\kappa}-E)\geq 0$ and this is valid for any compact $K$. But our compact set in consideration is  $K=\varphi (M)\cap (B_{N}(R)\times \mathbb{R})$ so that the boundary $\partial K\subset \partial B_{N}(R)\times \mathbb{R}$. Thus we have  $(\bar{E}_{\kappa}-E)\vert_{\partial K}= 0$ and then  $\bar{E}_{\kappa}\leq E$ in $K$.

\begin{rem}
\vspace{3mm}When $b<0$ the equality $\bar{E_{b}}=E$ in  $K$ implies that  $\vert \grad_{N} s\vert=\vert \grad_{K}s\vert$. Thus $\vert \grad_{N}\rho_{N}\vert =\vert \grad_{K}\rho_{N}\vert$ at every point of $K$. Recall that $\grad_{K}\rho_{N}=\sum_{i=1}^{m}\langle \grad_{N}\rho_{N},X_{i} \rangle X_{i}$.  Then  $\grad_{N} \rho_{N}$ is tangent to $K$. The  integral curves of $\grad_{N} \rho_{N}$ are geodesics (liftings of radial geodesics in $N$ via the projection map),  in $N\times \mathbb{R}$ and then in $K$. Thus we conclude that  through every point $q$ of $K$ passes a lifting   of a radial geodesic of $N$ passing through $\pi_{1}(q)$. Moreover, going through the computations (\ref{hessiano-de-s-A}) it is easy to see that   $$\sum_{i=2}^{m}\beta_{i}^{2}=1\,\; {\rm  and }\, \; (\gamma_{j}^{i})^{2}(\hess_{N}\rho_{N}(\partial/\partial \theta_{j}, \partial/\partial \theta_{j})-\displaystyle\frac{C_{b}}{S_{b}})=0 ,\,\; i=2,\ldots, m\,\;j=1, \ldots, n-1,$$ everywhere in $K$.
  
\end{rem}

\subsection{Proof of Theorem \protect\ref{thmA}} We start stating a lemma proved by Palmer. \begin{lem}[Palmer, \cite{palmer-1}]Let $E_{b}$ be the mean time exit of the ball $B_{\mathbb{N}^{m-1}(b)}(R)$. Then
\begin{equation} E_{b}'(t)=-\displaystyle\frac{{\rm vol}_{m-1}(B_{\mathbb{N}^{m-1}(b)}(t))}{{\rm vol}_{m-2}(\partial B_{\mathbb{N}^{m-1}(b)}(t))} \cdot \end{equation}
\end{lem}
\noindent Inequality (\ref{eq3.26}) says that $\triangle_{K}\bar{E}_{b}\circ \varphi \leq -1$ on $K$. Integration  over $K$ yields,
$$-{\rm vol}_{m}(K)=\int_{K}1\geq \int_{K}-\triangle_{K}\bar{E}_{b}\circ \varphi=-\int_{\partial K}\langle E_{b}'\, \grad \rho_{N}\circ \varphi, \nu\rangle\geq -\sup_{\partial K}\Vert E_{b}' \Vert {\rm vol}_{m-1}(\partial K) $$
Thus we have that
\begin{equation}\displaystyle \frac{{\rm vol}_{m-1}(\partial K)}{{\rm vol}_{m}(K)}\geq \frac{1}{\sup_{\partial K}\Vert E_{b}' \Vert}=\frac{{\rm vol}_{m-2}(\partial B_{\mathbb{N}^{m-1}(b)}(r_{K}))}{{\rm vol}_{m-1}(B_{\mathbb{N}^{m-1}(b)}(r_{K}))}\cdot \label{eq2.11}
\end{equation}

\subsubsection{Mean time exit on spherically symmetric manifolds.} 
A
 spherically symmetric manifold  is a quotient space $W=([0,R)\times
  \mathbb{S}^{n-1})/\backsim$,  $R\in (0,\infty]$, where $(t,\theta)\backsim (s,\alpha)$ iff $t=s $ and $\theta =\alpha$ or $s=t=0$,
  endowed with a Riemannian
metric of the form
 $dt^{2}+f^{2}(t)d\theta^{2}$, where $f\in C^{2}([0,R])$  with
  $f(0)=0$, $f'(0)=1$, $f(t)>0$ for
all $t\in (0,R]$. 
 The class of spherically
symmetric manifolds includes the canonical space forms
$\mathbb{R}^{n}$, $\mathbb{S}^{n}(1)$ and $\mathbb{H}^{n}(-1)$. Let $B_{W}(r)\subset W$ be a geodesic ball of radius $r$ and center $0=(0\times \mathbb{S}^{n-1})/\backsim $ in a spherically  symmetric manifold $(W, dt^{2}+f^{2}(t)d\theta^{2})$. The mean time of the first exit $E$ of $B_{W}(r)$ is given by
\begin{equation}E(x)=E(\vert x\vert)=-\int_{\vert x\vert}^{r}\frac{1}{f^{n-1}(\sigma)}\int_{0}^{\sigma}f^{n-1}(s) dsd\sigma. \end{equation} as one can easily check that $E$ satisfies $\triangle_{W}E=-1$ in $B_{W}(r)$ with $E\vert \partial B_{W}(r)=0$. Here $\vert x\vert = {\rm dist}_{W}(0, x)$. It is also straightforward to show that $$E'(\vert x\vert)=-\frac{\int_{0}^{r}f^{n-1}(s) ds}{f^{n-1}(r)}=-\frac{\vol_{n}(B_{W}(r))}{\vol_{n-1}(\partial B_{W}(r))}\cdot$$  Consider the Dirichlet problem $ \triangle_{W}u+\lambda_{1}(B_{W}(r))u=0$ in $B_{W}(r)$ with $u=0$ on $\partial B_{W}(r)$.
 It was shown by the authors in \cite{bessa-montenegro0} that the first Dirichlet eigenvalue $\lambda_{1}(B_{W}(r))$ is bounded below by \begin{equation} \lambda_{1}(B_{W}(r))\geq \frac{[\inf_{B_{W}(r)}\diver X]^{2}}{4\sup_{B_{W}(r)}\vert X\vert^{2}},\label{eqLambda}\end{equation}where $X$ is a vector field in $B_{W}(r)$ with $\inf \diver X >0$ and $\sup \vert X\vert <\infty$. Taking $X=-\grad_{W}E$ we have that $\diver X=1$ and $\vert X\vert =\vert E'\vert$. Applying (\ref{eqLambda}) we obtain the following theorem.
 \begin{thm} Let $B_{W}(r)$ be a geodesic  ball centered at  $0=(0\times \mathbb{S}^{n-1})/\backsim $ with of radius $r$ in a spherically symmetric Riemannian $n$-manifold $(W, dt^{2}+f^{2}(t)d\theta^{2})$. Let $V(t)$ and $S(t)$ be respectively the $n$-volume and $(n-1)$-volume of $B_{W}(t)$ and $\partial B_{W}(t)$. Then \begin{equation}\label{eqLambda2}\lambda_{1}(B_{W}(r))\geq \inf_{0\leq t\leq r}\frac{1}{4} \left[\frac{S(t)}{V(t)}\right]^{2}
 \end{equation}
 \end{thm}
 \begin{cor}Let $W$ be a complete non-compact spherically symmetric manifold. Suppose that the boundary of $B_{W}(t)$ has  volume growth $ c_{1}e^{c_{3}t}\leq S(t)\leq c_{2}e^{c_{3}t}$, $c_{1}<c_{2}$ and $c_{3}$ are positive constants. Then $$\lambda^{\ast}(W)=\lim_{r\to \infty}\lambda_{1}(B_{W}(r))\geq ( \frac{c_{1}c_{3}}{2c_{2}})^{2}\cdot$$
 \end{cor}
\begin{rem}The inequality (\ref{eqLambda2}) should be compared with the inequality $\lambda_{1}(B_{W}(r))\geq \displaystyle \frac{1}{\int_{0}^{r}\frac{V(\sigma )}{S(\sigma)}d\sigma}$ proved by  Barroso and Bessa  in \cite{barroso-bessa}.
\end{rem}

\subsection{More on submanifolds with tamed second fundamental form}  If $\varphi :M\hookrightarrow N$ is a complete $m$-submanifold with  tamed second fundamental form  immersed in a  Hadamard $n$-manifold  with sectional curvature $K_{N}\leq b\leq 0$ then $\varphi $ is proper. Moreover $\varphi (M)$ has finite topology, see \cite{bessa-jorge-montenegro}, \cite{costa}. In this section we are going to present the idea to prove that $\varphi (M)$ has finite topology since we need a corollary from its proof. Recall that a submanifold  $\varphi :M\hookrightarrow N$ has  tamed second fundamental form if $\lim_{i\to \infty}a_{i}(M)=a(M)<1$. Thus given $c\in (a(M),1)$ there is an $r_{o}>0$ such that $$\Vert \alpha(\varphi (x)) \Vert \leq c\cdot \frac{C_{b}}{S_{b}}(\rho_{N}(\varphi(x))$$ for all $x\in M\setminus B_{M}(r_{0})$. Here $\rho_{N}$ is the intrinsic distance function in $N$ to a point $p=\varphi (q)$.
Let $r>r_{0}$ be such that $\varphi (M)\pitchfork \partial B_{N}(r)$ and let $\Gamma = \varphi (M)\cap \partial B_{N}(r)$. Setting $\Lambda=\varphi ^{-1}(\Gamma)$ we construct a smooth vector field $\nu $ on a open neighborhood of $\Lambda$  so that $\forall x\in \Lambda$, $y=\varphi (x)$ we have that $$T_{y}M=T_{y}\Gamma \oplus [[d\varphi (x).v(x), \,\grad \rho_{N}]]$$ with $\langle d\varphi (x).\nu(x), \grad \rho_{N}\rangle >0$. Here $[[d\varphi (x).\nu(x), \,\grad \rho_{N}]]$ is the vector space generated by $d\varphi (x).\nu(x)$ and $\grad \rho_{N}$. For simplicity of notation we are going to identify   $d\varphi (x).\nu(x)=\nu (y)$. Define $\psi (x)=\langle \nu (y), \grad_{N}(y) \rangle$, $x\in \Lambda$. Since $\Lambda$ is compact and $\psi (x)>0$ there is a positive minimum $\psi_{0}$. Consider the Cauchy Problem on $M$ \begin{equation}\label{cauchy}
\left\{
\begin{array}{ccl}
\xi_t(t,x)&=&\displaystyle\frac{1}{\psi}\,\nu(\xi(t,x))\\
\\
\xi(0,x) &=& x\\
\end{array}
\right.
\end{equation} It was shown in \cite{costa} that $\psi $ satisfies the following differential equation along the integral curves $\xi (t,x)$
\begin{equation}\label{eqdif1}
\begin{array}{ccl}
-(\sqrt{1 - \psi^2})_t &=&\sqrt{1 - \psi^2} \; \hess\rho_{N}(\omega,\omega) + \langle \nu^*,\alpha(\nu,\nu)\rangle\\
\end{array}
\end{equation}Where $\nu^*$ is a unit vector normal to $\nu$ and $\omega$ is a unit normal vector to $TM$ and to $\grad \rho_{N}$.
As consequence of Hessian Comparison Theorem we have the  following inequality
\begin{equation}\label{eqdif2}
-(\sqrt{1 - \psi^2})_t \geq \sqrt{1 - \psi^2}\; \displaystyle \frac{C_{b}}{S_{b}}+ \langle
\nu^*,\alpha(\nu,\nu)\rangle\\
\end{equation}  It was shown that $\rho_{N}(\xi (t, x))=t+r_{0}$ so  we can write $S_{b}(\rho_{N}(\xi(t,x)))=S_{b}(t+r_{0})$. The inequality (\ref{eqdif2}) is equivalent to
\begin{equation}\label{eqdif4}
\left[S_{b}(t + r_0)\sqrt{1 - \psi^2}\right]_t \leq- S_{b}(t + r_0) \langle \nu^*,\alpha(\nu,\nu)\rangle
\end{equation}Integrating  $(\ref{eqdif4})$ of $0$ to $t$ we obtain
$$ \sqrt{1-\psi^{2}}(t)\leq \frac{S_{b}(r_{0})}{S_{b}(t+r_{0})} \sqrt{1-\psi^{2}}(0) -\frac{1}{S_{b}(t+r_{0})}\int_{0}^{t} S_{b}(t+r_{0})\langle
\nu^*,\alpha(\nu,\nu)\rangle ds
$$But $-\langle
\nu^*,\alpha(\nu,\nu)\rangle(\xi(s,x))\leq \Vert \alpha (\xi(s,x)) \Vert \leq c\cdot (C_{b}/S_{b})(s+r_{0})$. Thus

\begin{equation} \label{eqdif5}\sqrt{1-\psi^{2}}(t)\leq  \frac{S_{b}(r_{0})}{S_{b}(t+r_{0})} \left[\sqrt{1-\psi_{0}^{2}}-c\right] +c<1
\end{equation} Where $\psi_{0}= \inf_{\Lambda} \psi >0$. Thus for every $x\in \Lambda$ we have that the function $\psi$ satisfies 
inequality (\ref{eqdif5}) along the integral curve $\xi(t,x)$, (which is defined for all $t$).

  \vspace{3mm}

  \noindent Now let $f:M\to \mathbb{R}$ be defined by $\rho_{N}\circ \varphi$. The gradient of $ f$ is the projection of $\grad \rho_{N}$ on $TM$, i.e.  $\grad f = \langle \grad \rho_{N}, \nu\rangle \nu=\psi \cdot \nu$. Set
 \begin{equation} \label{eqB}B(b,c,r_{0})=\sup_{t,x}\displaystyle \frac{S_{b}(r_{0})}{S_{b}(t+r_{0})} \left[\sqrt{1-\psi_{0}^{2}}-c\right] +c<1\end{equation} By (\ref{eqdif5}) we have that $$\inf_{M\setminus \varphi^{-1}(B_{N}(r_{0}))} \psi  \geq \sqrt{1-B^{2}(b,c,r_{0})}>0  $$ Therefore we have that $\Vert \grad f (x)\Vert \geq \sqrt{1-B^{2}(b,c,r_{0})}$, for $x\in M\setminus \varphi^{-1}(B_{N}(r_{0}))$.

\subsection{Proof of Theorem \ref{thm2}.}Let $\varphi :M\hookrightarrow N$ be a complete $m$-dimensional submanifold of a complete Hadamard manifold with sectional curvature $b_{1}\leq K_{N}\leq b_{2}\leq 0$ with tamed second fundamental form . As we mentioned  before $\varphi $ is proper. Fix a point $p=\varphi (q)\in N$ and $c\in (a(M),1)$ let $D(R)$ be the extrinsic geodesic ball with center at $p$ and radius $R\geq r_{0}=r_{0}(c)>0$. Consider $f:M\to \mathbb{R}$ given by $f=\rho_{N}^{2}\circ \varphi$. We want to estimate $\sup_{D(R)}\triangle f$. We proceed as follows.
 \begin{eqnarray}\triangle f &=&  2\sum_{i}^{m}\langle X_{i},\, \grad \rho_{N}\rangle^{2} +2\rho_{N} \hess \rho_{N}(X_{i}, X_{i}) + 2 \rho_{N}\langle \grad \rho_{N}, \stackrel{\to}{H}\rangle\nonumber \\
 &\leq & 2\left[1+ \sup_{t\in [0,R]}\rho_{N} \frac{C_{b_{1}}}{S_{b_{1}}}(\rho_{N})\right]+\max\{2r_{0}\sup_{B_{N}(r_{0})}\vert H\vert,\, 2c\cdot R\cdot (C_{b_{2}}/S_{b_{2}})(R)\}\nonumber\\
 &=&2\left[1+ \sup_{t\in [0,R]}\rho_{N} \frac{C_{b_{1}}}{S_{b_{1}}}(\rho_{N})\right] + \Lambda (\sup_{B_{N}(r_{0})}\vert H\vert, c, r_{0}, R, b_{2})
 \end{eqnarray}

 \vspace{3mm} \noindent On the other hand $\grad f= 2\rho_{N} \psi  \nu$. Thus we have that $$\inf_{\partial D(R)}\Vert \grad f\Vert \geq 2\cdot R\cdot \sqrt{1-B^{2}(b_{2}, c,r_{0})}  $$  By Green's Theorem we have that
 \begin{eqnarray}\sup_{D(R)}\triangle f \cdot {\rm vol}_{m}(D(R))&\geq& \int_{D(R)} \triangle f
  =  \int_{\partial D(R)}\langle \grad f, \nu\rangle\nonumber \\ && \nonumber \\ &\geq & \inf_{\partial D(R)}\Vert \grad f\Vert \cdot {\rm vol}_{m-1}(\partial D(R))\nonumber \end{eqnarray}

 \vspace{2mm}

 From we obtain
 $$\frac{{\rm vol}_{m-1}(\partial D(R))}{{\rm vol}_{m}(D(R))}\leq \frac{1+ \sqrt{-b_{1}}\cdot R\cdot \coth(\sqrt{-b_{1}}\cdot R)+\Lambda (\sup_{B_{N}(r_{0})}\vert H\vert,c, r_{0}, R, b_{2})}{R\cdot \sqrt{1-B^{2}(b_{2}, c,r_{0})}}
$$


\subsection*{Acknowledgment}
We thank the Abdus Salam International Center for Theoretical Physics-ICTP for their generous hospitality and support. We were also partially supported by  CNPq-FUNCAP.
\end{document}